# NUMERICAL MODELLING OF VARIABLE DENSITY SHALLOW WATER FLOWS WITH FRICTION TERM

Hanini, A[1,4], Beljadid, A[1,2,5], Ouazar, D[3,6]
[1] Mohammed VI Polytechnic University, Morocco
[2] University of Ottawa, Canada
[3] Mohammadia School of Engineering, Morocco
[4] amine.hanini@um6p.ma
[5] abeljadi@uottawa.ca
[6] ouazar@emi.ac.ma

**Abstract:** In this study, we focus on the modelling of coupled systems of shallow water flows and solute transport with source terms due to variable topography and friction effect. Our aim is to propose efficient and accurate numerical techniques for modelling these systems using unstructured triangular grids. We used a Riemann-solver free method for the hyperbolic shallow water system and a suitable discretization technique for the bottom topography. The friction source term is discretized using the techniques proposed by (Xia and Liang 2018). Our approach performs very well for stationary flow in the presence of variable topography, and it is well-balanced for the concentration in the presence of wet and dry zones. In our techniques, we used linear piecewise reconstructions for the variables of the coupled system. The proposed method is well-balanced, and we prove that it exactly preserves the nontrivial steady-state solutions of the coupled system. Numerical experiments are carried out to validate the performance and robustness of the proposed numerical method. Our numerical results show that the method is stable, well-balanced and accurate to model the coupled systems of shallow water flows and solute transport.

## 1. INTRODUCTION

Numerical modelling of the scalar transport process is of a great importance for surface water pollution risk assessment which is used in a large variety of applications for environmental protection and management (Rekolainen et al. 2003, Horn et al. 2004). The coupled system of shallow water equations and solute transport model are used to predict the dynamics of both water flow and contaminant transport in rivers and coastal regions (Cea and Vazquez-Cendon 2012, Begnudelli and Sanders 2006). In shallow flows, horizontal scales are predominant over vertical ones and the pressure is nearly hydrostatic (Vreugdenhil 1994). The shallow water equations can be derived by integrating vertically the three-dimensional incompressible Navier-Stokes equations and assuming hydrostatic pressure distribution where vertical acceleration and viscous effects are neglected.

The design of efficient numerical schemes for solving the coupled model of shallow water flow and scalar transport remains major challenge which requires robust numerical techniques (Murillo et al. 2012) to ensure important physical properties of schemes such as the well-balanced property (Audusse et al. 2004). Among the difficulties that may be encountered, the resulting nonlinear hyperbolic system can lead to non-smooth solutions which can exhibit singularities, and in the presence of source terms, computational techniques may lead to numerical oscillations caused by the imbalance between source and flux terms. For instance, the friction source term which is non-linear can leads to numerical issues (Xia et al 2017). The important challenge here for the proposed numerical method is to discretize properly this term. In our approach, we used the implicit techniques proposed by (Xia and Liang 2018) to discretize the friction term. These techniques can be implemented explicitly, and our aim is to efficiently solve the coupled system of shallow water flow and solute transport with friction term.

Many advantages allow the finite volume method to be convenient for solving system of conservation or balance laws, for instance, this method is conservative (Audusse et al. 2004, Beljadid et al. 2012, Bryson et al. 2011), and deals with solutions with large gradients. One of the most successful numerical approach to solve hyperbolic PDE in the finite volume framework is due to Godunov (Godunov 1959) and is based on representing the solution at cell centers on both sides of each cell interface of the whole mesh by piecewise constant states. To have more accurate results, instead of using constant states, we may consider polynomial ones for the reconstruction variable in each cell. For instance, linear piecewise approximation (Van Leer 1979) can be used to design second order accurate schemes. Thereafter, many methods have been developed such as central-upwind schemes (Kurganov 2018, Kurganov and Levy 2002) which are based on characteristic information on local speeds of propagation on the two sides of cell interfaces. These numerical schemes have been used in many studies to solve the shallow water system and have the advantage of simplicity of implementation and lead to accurate results for this system.

The outline of this paper is as follows. In Section 2, we present the coupled system of shallow water equations and scalar transport model and the semi-discrete form of the numerical scheme used to solve this system. In Section 3 the well-balanced discretization of the source term due to variable bottom topography is discussed as well as the discretization of friction source term. We present numerical experiments in Section 4 to test the performance of the proposed scheme. Finally, we provide some concluding remarks in Section 5.

## 2. GOVERNING EQUATIONS

In this study, we used the following coupled system of shallow water equations (SWEs) with variable density and transport model. The SWEs, also called the Saint-Venant equations can be derived from the depth-averaged three-dimensional incompressible Navier-Stokes equations (Vreugdenhil 1994). To model the solute transport in water, we used the depth-averaged scalar transport equation.

$$[1] \begin{cases} \partial_t hr + \partial_x hur + \partial_y hvr = 0 \\ \partial_t hur + \partial_x \left(hu^2 r + \frac{1}{2}gh^2 r\right) + \partial_y huvr = -ghr\,\partial_x Z - gn^2 (h)^{-1/3} u\sqrt{(u^2 + v^2)} \\ \partial_t hvr + \partial_y huvr + \partial_x \left(hv^2 r + \frac{1}{2}gh^2 r\right) = -ghr\,\partial_y Z - gn^2 (h)^{-1/3} v\sqrt{(u^2 + v^2)} \\ \partial_t hc + \partial_x huc + \partial_y hvc = 0, \end{cases}$$

where $x$ and $y$ are the spatial coordinates, $t$ is time, $h$ is the mixture depth, $u$ and $v$ are the depth-averaged velocities of the flow in the $x-$ and the $y-$directions, $c$ is the scalar depth-averaged volumetric concentration, $Z$ is the elevation of the bottom topography, $g$ is the acceleration due to gravity. $r$ is the relative density of the mixture to that of the clean water and $n$ is the Manning coefficient (Manning 1891).

Equation 1 can be rewritten into the following vectorial form

$$[2]\ \partial_t \mathbf{U} + \partial_x \mathbf{F}(\mathbf{U}) + \partial_y \mathbf{G}(\mathbf{U}) = \mathbf{S}_b + \mathbf{S}_f,$$

where we use the vector of variables $\mathbf{U} = (q_1, q_2, q_3, q_4)^T$, $q_1 \coloneqq hr$, $q_2 \coloneqq hur$, $q_3 \coloneqq hvr$, and $q_4 \coloneqq hc$. We denote by $w \coloneqq h + Z$ the water surface elevation. The friction source term $\mathbf{S}_f$ is expressed using the manning formula and $\mathbf{S}_b$ is the source term due to variable bottom topography.

The fluxes in the x- and y-directions are as follows:

$$\mathbf{F}(\mathbf{U}) = \left(q_2, \frac{q_2^2}{q_1} + \frac{1}{2}gq_1(q_1 - \Delta q_4), \frac{q_2 q_3}{q_1}, \frac{q_2 q_4}{q_1}\right)^T$$

$$\mathbf{G}(\mathbf{U}) = \left(q_3, \frac{q_2 q_3}{q_1}, \frac{q_3^2}{q_1} + \frac{1}{2}gq_1(q_1 - \Delta q_4), \frac{q_3 q_4}{q_1}\right)^T.$$

The bottom topography and friction source terms are, respectively, expressed as

$$\mathbf{S}_b = \left(0, -g(w-Z)r\, \partial_x Z, -g(w-Z)r\, \partial_y Z, 0\right)^T$$

$$\mathbf{S}_f = \left(0, -gn^2(w-Z)^{-\frac{7}{3}}\frac{q_2}{r^2}\sqrt{q_2^2 + q_3^2}, -gn^2(w-Z)^{-\frac{7}{3}}\frac{q_3}{r^2}\sqrt{q_2^2 + q_3^2}, 0\right)^T,$$

We note that in the case of small water depth $h \approx 0$, the friction source term $\mathbf{S}_f$ becomes stiff which is very challenging in terms of numerical discretization.

The relative density of the mixture to that of the clean water is related to the volumetric concentration via the following formula

[3] $r = 1 + \Delta c,$

where $\Delta = \frac{\rho - \rho_w}{\rho_w}$ is the relative density of the constituent, $\rho$ is its density, $\rho_w r$ is the mixture density and $\rho_w$ is the density of water.

We used a triangulation $\mathcal{T} := \bigcup_j T_j$ of the computational domain, using non-overlapping triangular cells $T_j$ of measure $|T\_j|$. We denote by $(x_j, y_j)$ the coordinates of the center of mass $G_j$ of triangle $T_j$, $M_{jk} = (x_{jk}, y_{jk})$ is the midpoint of the k-th side of the triangle, $T_{j1}$, $T_{j2}$ and $T_{j3}$ are the three neighboring triangles sharing a common side with $T_j$, and $n_{jk} := (\cos(\theta_{jk}), \sin(\theta_{jk}))^T$ is the outer unit normal vector to the k-th side $\Gamma_{jk}$ of $T_j$, of length $|\Gamma_{jk}|$, with k = 1,2,3 as shown in Figure 1.

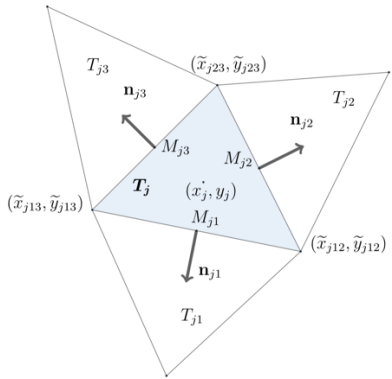

Figure 1: A typical triangular control volume and three neighbors

The variable $\bar{\mathbf{U}}_j(t)$ represents the approximation of the cell average of the solution at time $t$

$$\overline{U}_j(t) \approx \frac{1}{|T_j|}\int_{|T_j|} U(t,x,y)\,dxdy.$$

We consider the semi-discrete form of numerical scheme [4], which is described by the following ODE

$$[4]\quad \frac{d\overline{U}_j}{dt} = -\frac{1}{|T_j|}\sum_{k=1}^{3}\frac{|\Gamma_{jk}|\cos(\theta_{jk})}{a_{jk}^{in}+a_{jk}^{out}}\left[a_{jk}^{in}F\left(U_{jk}(M_{jk})\right)+a_{jk}^{out}F\left(U_j(M_{jk})\right)\right]$$

$$-\frac{1}{|T_j|}\sum_{k=1}^{3}\frac{|\Gamma_{jk}|\sin(\theta_{jk})}{a_{jk}^{in}+a_{jk}^{out}}\left[a_{jk}^{in}G\left(U_{jk}(M_{jk})\right)+a_{jk}^{out}G\left(U_j(M_{jk})\right)\right]$$

$$+\frac{1}{|T_j|}\sum_{k=1}^{3}|\Gamma_{jk}|\frac{a_{jk}^{in}a_{jk}^{out}}{a_{jk}^{in}+a_{jk}^{out}}\left[U_{jk}(M_{jk})-U_j(M_{jk})\right]+(\overline{S}_b)_j+(\overline{S}_f)_j.$$

The values $U_j(M_{jk})$ and $U_{jk}(M_{jk})$ of the variables of the system [1] at the midpoints of cell interfaces, respectively at the left and right side of each interface of the computational domain are approximated using the linear piecewise reconstructions $\widetilde{U}$

$$U_j(M_{jk}) = \lim_{(x,y)\to M_{jk},(x,y)\in T_j}\widetilde{U}(x,y),$$

$$U_{jk}(M_{jk}) = \lim_{(x,y)\to M_{jk},(x,y)\in T_{jk}}\widetilde{U}(x,y),$$

In this study, we used a linear reconstruction of the discrete variables of the system [4] defined by

$$[5]\quad \widetilde{U}(x,y) := \overline{U}_j + (\nabla_x U)_j(x-x_j)+(\nabla_y U)_j(y-y_j). \qquad (x,y)\in T_j.$$

In Equation 5, $\overline{U}_j$ corresponds to the value $U_j := \widetilde{U}(x_j,y_j)$ at the center of mass of the cell $T_j$.

The directional local speeds at interfaces between two cells are approximated based on the eigenvalues of the Jacobian matrix $J_{jk}$ of the system:

$$\begin{cases} a_{jk}^{in} = -\min\{\lambda_{min}\left[J_{jk}\left(U_j(M_{jk})\right)\right],\lambda_{min}\left[J_{jk}\left(U_{jk}(M_{jk})\right)\right],0\} \\ a_{jk}^{out} = \max\{\lambda_{max}\left[J_{jk}\left(U_j(M_{jk})\right)\right],\lambda_{max}\left[J_{jk}\left(U_{jk}(M_{jk})\right)\right],0\} \end{cases}$$

where $\lambda_{min}$ and $\lambda_{max}$ are respectively the minimum and maximum eigenvalues of $J_{jk}$ and

$$J_{jk} = \frac{\partial F}{\partial U}\cos(\theta_{jk})+\frac{\partial G}{\partial U}\sin(\theta_{jk}).$$

The expressions of the eigenvalues and eigenvectors of the matrix $J_{jk}$ are given in (Hanini et al. 2021). For simplicity, the following normal velocities at the midpoints $M_{jk}$ are introduced

$$u_j^\theta(M_{jk}) = \cos(\theta_{jk})u_j(M_{jk})+\sin(\theta_{jk})v_j(M_{jk}), \qquad u_{jk}^\theta(M_{jk}) = \cos(\theta_{jk})u_{jk}(M_{jk})+\sin(\theta_{jk})v_{jk}(M_{jk}),$$

which are used to express the directional local speeds as follows

$$\begin{cases} a_{jk}^{in} = -\min\{u_j^\theta(M_{jk})-\tilde{c}_j(M_{jk}), u_{jk}^\theta(M_{jk})-\tilde{c}_{jk}(M_{jk}),0\} \\ a_{jk}^{in} = -\min\{u_j^\theta(M_{jk})+\tilde{c}_j(M_{jk}), u_{jk}^\theta(M_{jk})+\tilde{c}_{jk}(M_{jk}),0\} \end{cases}$$

where $\tilde{c}_j(M_{jk}) = \sqrt{gh_j(M_{jk})}$ and $\tilde{c}_{jk}(M_{jk}) = \sqrt{gh_{jk}(M_{jk})}$.

To improve the accuracy of the numerical method, we used a linear reconstruction for the variables of the system where the gradient of the ith component of **U** is computed using Green's formula

[6] $\nabla \mathbf{U}^i \approx \frac{1}{|T_j|} \int_{T_j} \nabla \mathbf{U}^i dxdy = \frac{1}{|T_j|} \int_{\partial T_j} \mathbf{U}^i \mathbf{n} d\sigma \approx 1|T_j| \sum_{k=1}^{3} |\Gamma_{jk}| \mathbf{U}^i(M_{jk}) \mathbf{n}_{jk}$,

where **n** is the unit normal vector to the boundary $\partial T_j$.

In Equation 6, we choose $\mathbf{U}^i(M_{jk}) = \frac{1}{2}(\mathbf{U}^i_j + \mathbf{U}^i_{jk})$ to obtain the approximation of the gradients of the variables in each cell of the computational domain. To avoid oscillations in numerical solutions, we introduce a minmod limiter function (Nessyahu and Tadmor 1990) which is applied to the components of the computed cell gradients. In the $x$-direction, we use the following correction:

$$\frac{\partial^{lim}\mathbf{U}_j}{\partial x} = \frac{1}{2}\left(\min_{k\in\mathcal{N}(j)} \text{sign}\left[\frac{\partial \mathbf{U}_k}{\partial x}\right] + \max_{k\in\mathcal{N}(j)} \text{sign}\left[\frac{\partial \mathbf{U}_k}{\partial x}\right]\right) \min_{k\in\mathcal{N}(j)} \text{sign}\left|\frac{\partial \mathbf{U}_k}{\partial x}\right|,$$

where $sign$ is the standard sign function, $\mathcal{N}(j)$ is the set of neighboring cells sharing common side with the cell $T_j$. Similar steps are used to obtain the component of the cell gradients in the $y$-direction $\frac{\partial^{lim}\mathbf{U}_j}{\partial x}$.

The new approximation of the cell gradient is used in Equation 5.

## 3. DISCRETIZATION OF THE TOPOGRAPHY AND FRICTION SOURCES TERMS

The cell average of the source term due to bottom topography $\overline{\mathbf{S}}_j \approx \frac{1}{|T_j|} \int_{T_j} \mathbf{S}(t,x,y) dxdy$ is discretized using the approach developed in (Hanini et al. 2021) where variable density is taken into account. The semi-discrete form of the numerical scheme [4] is well-balanced and it exactly preserves the following steady-states solutions

[7] $\quad h \equiv h_0 = \left(\frac{n^2 q_0^2}{r_0^3 C}\right)^{3/10} \quad q_2 \equiv q_0 \quad q_3 \equiv 0 \quad r \equiv r_0 \quad \partial_x Z = -C \quad \partial_y Z = 0 \quad c \equiv c\_0$,

where $q_0$, $r_0$, $C(C > 0)$ and $c_0$ are constants. In the following, we will prove that [7] is an exact solution of the ODE [4].

We note that since we are dealing with stationary solutions, the term in the left-hand side of Equation 4 is zero. From Equation 7 $q_3 \equiv 0$, then all the components of the flux **G** introduced in Section. 2 vanish and consequently the second term of the right-hand side of Equation 4 is zero. For the steady state solution [7] we have $\mathbf{U}_{jk}(M_{jk}) = \mathbf{U}_j(M_{jk})$, then the third term on the right-hand side of Equation 4 vanishes. Finally, Equation 4 becomes

$$0 = -\frac{1}{|T_j|} \sum_{k=1}^{3} \frac{|\Gamma_{jk}| \cos(\theta_{jk})}{a_{jk}^{in} + a_{jk}^{out}} \left[(a_{jk}^{in} + a_{jk}^{out}) \mathbf{F}\left(\mathbf{U}_j(M_{jk})\right)\right] + (\overline{\mathbf{S}}_b)_j + (\overline{\mathbf{S}}_f)_j,$$

which gives the following formulation

$$0 = -\frac{1}{|T_j|}\sum_{k=1}^{3}|\Gamma_{jk}|\cos(\theta_{jk})\left[\frac{q_0^2}{q_1}+\frac{1}{2}gh_0^2 r_0\right] + \frac{g}{2|T_j|}r_0\sum_{k=1}^{3}|\Gamma_{jk}|\cos(\theta_{jk})h_0^2 + gr_0 h_0 C - gn^2(h_0)^{-\frac{7}{3}}\frac{q_0}{r_0^2}|q_0|$$

then by simplification, we obtain

$$0 = -\frac{1}{|T_j|}\sum_{k=1}^{3}|\Gamma_{jk}|\cos(\theta_{jk})\left[\frac{q_2^2}{q_1}\right] + gr_0 h_0 C - gn^2(h_0)^{-\frac{7}{3}}\frac{q_0}{r_0^2}|q_0|$$

This leads to

$$0 = -\frac{1}{|T_j|}\sum_{k=1}^{3}|\Gamma_{jk}|\cos(\theta_{jk})\left[\frac{q_2^2}{q_1}\right] = -\frac{1}{|T_j|}\left[\frac{q_2^2}{q_1}\right]\sum_{k=1}^{3}|\Gamma_{jk}|\cos(\theta_{jk})$$

This equation is valid since in every triangle the relation $\sum_{k=1}^{3}|\Gamma_{jk}|\cos(\theta_{jk}) = 0$ is verified, which can be obtained using divergence theorem.

Following the same steps for the second momentum equation in Equation 4, one obtains

$$0 = -\frac{1}{|T_j|}\sum_{k=1}^{3}|\Gamma_{jk}|\sin(\theta_{jk})\left[\frac{q_3^2}{q_1}\right].$$

The friction source term $S_f$ is discretized using the approach developed in (Xia and Liang 2018) which allows us to obtain a semi-implicit scheme. The obtained system is solved explicitly to avoid the use of computationally expensive iterative methods.

## 4. NUMERICAL EXAMPLES

Here, we perform two numerical examples to test the performance of the numerical scheme. A first-order forward Euler method for temporal discretization is adopted with the following stability CFL condition

$$\Delta t \leq \frac{1}{6}\frac{\min_{j,k}\{H_{jk}\}}{\max_{j,k}\{a_{jk}^{out}\}},$$

where $H_{jk}$ are the heights of the triangle $T_j$. In all numerical tests the constant of gravitation is $g = 9.81.\,ms^{-2}$ and CFL= 0.8.
In the case when the water depth is very small or even zero, to calculate the velocities $u$, $v$ and depth averaged concentration $c$ we used the desingularisation methodology detailed in (Kurganov and Petrova 2007).

### 4.1 Example. 1 : Steady flow over a slanted surface

This example has been adopted to verify the ability of our numerical scheme to preserve the steady state solutions [7] in the presence of variable topography as shown in Figure. 2. We consider a computational domain $[0, 10] \times [0, 5]$ discretized with a non-uniform triangular cells with a cell average area $2.\,10^{-3}$ subject to the constant initial data [7] with $h_0 = 0.25119\,m$, $q_0 = 0.1\,m^2.s^{-1}$, $n = 0.1\,m^{-1/3}.s$, $C = -0.01$, $\Delta = 0.1$ wich gives $r = 1.1$. Wall boundary conditions are set in the $y-$direction while outflow type-condition are set upstream and downstream of the domain. The numerical simulations are performed until $t = 100\,s$.

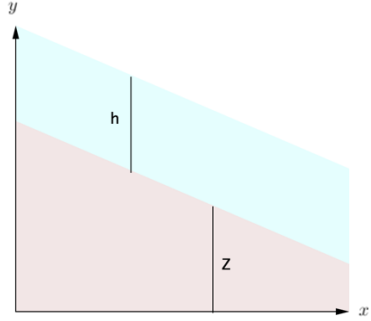

Figure 2: Bottom and water surface level for example 1

Table 1: Errors in computing the steady states [7] for example 1

|   | $L_1$ −error | $L_2$ −error | $L_\infty$ −error |
|---|---|---|---|
| $h$ | $1.99 \times 10^{-4}$ | $2.29 \times 10^{-4}$ | $2.42 \times 10^{-4}$ |
| $q_0$ | $8.37 \times 10^{-4}$ | $1.1 \times 10^{-3}$ | $1.4 \times 10^{-3}$ |
| $c$ | $1.3 \times 10^{-3}$ | $1.3 \times 10^{-3}$ | $3.9 \times 10^{-3}$ |

Table 1 shows the computed errors and Figure. 3 presents the computed solution for the water surface level $w$, $x$−discharge $q_2$, concentration $c$ and $x$−velocity $u$ at time $t = 100\ s$. The results show that the computed solution is free of any oscillation. The numerical solution remains stationary which confirm the robustness of the proposed numerical scheme.

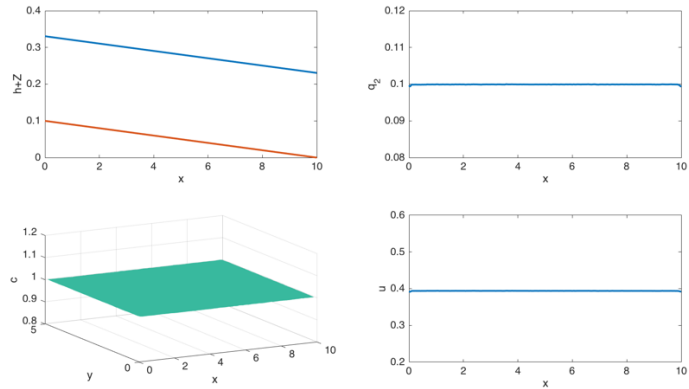

Figure 3: Solution $(w, Z, q_2, c, u)$ computed at time $t = 100\ s$

## 4.2  Example. 2 : Partial variable density dam-break over dry bed

In this numerical example, we consider a 2D partial dam-break with a symmetrical breach. We use the computational domain $500 \times 300$ as shown in Figure 4. This numerical test is performed using wall-boundary condition at all sides of the computational domain which is discretized using triangular meshes with an average cell area $|Tj| = 4.68$. Initially the upstream water depth is 10 and volumetric concentration is 1, and they are both equal to 0 downstream. The relative density is $\Delta = 0.1$ which gives a density of the mixture of 1.1. The Manning roughness coefficient is set to 0.1. The snapshots of the results displaying the water height and volumetric concentration at time $t = 10000\ s$ after dam failure. The numerical solution is free of any oscillation and remains symmetric. The scheme remains stable despite the stiffness of the friction source term for small values of the water depth. The concentration remains constant in the wet and dry region where $c = 1$ for $h > 0$. and $c = 0$ for $h = 0$ m. At the final simulation time, the concentration is constant $c = 1$ which confirms that the numerical scheme is well-balanced in terms of concentration.

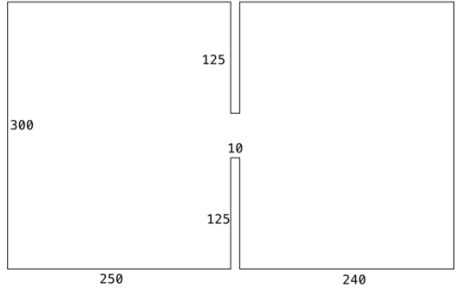

Figure 4: Sketch of the domain for example 2

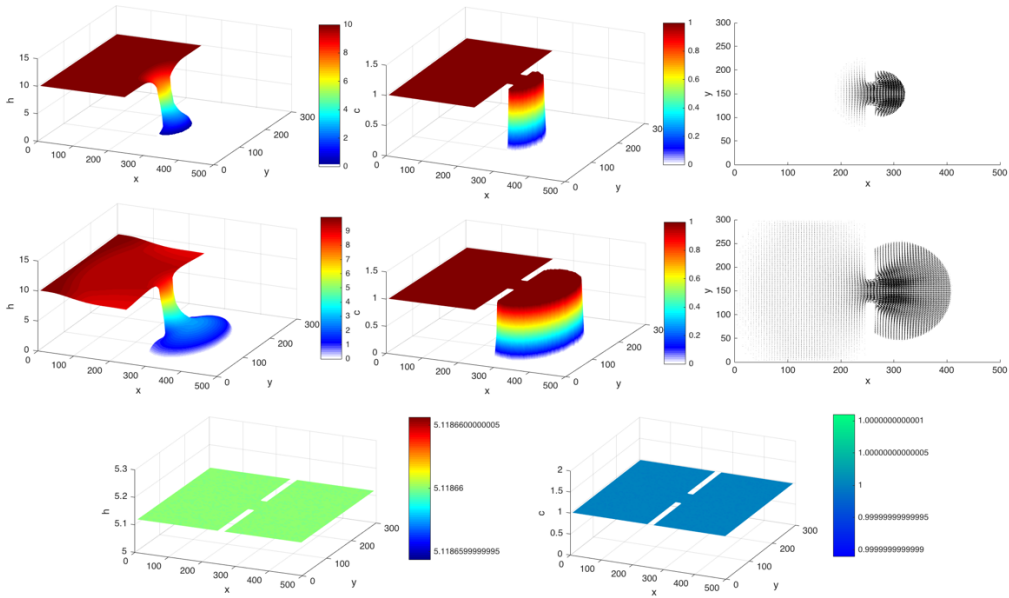

Figure 5: Snapshots of the solution of the water height $h$, volumetric concentration $c$ and velocity field at $t = 5\ s$ (first row), $t = 25\ s$ (second row) and $t = 10000\ s$ (long-time simulation, third row).

## 5.  CONLUSIONS

In this paper, we have studied a coupled 2D shallow water system and solute transport model with variable topography and friction source term. We proposed accurate numerical techniques for modelling this coupled system using triangular grids. A central-upwind scheme is used to solve the hyperbolic shallow water equation. A well-balanced discretization technique is used for the bottom topography. For accurately discretize the friction source term, we applied the techniques proposed by (Xia and Liang 2018). Our approach performs very well for nontrivial stationary flow in the presence of variable topography. The proposed scheme is well-balanced, and we prove that it exactly preserves the nontrivial steady-state solutions. Our numerical experiments demonstrate the robustness of the propose techniques and confirm that the proposed method is stable, well-balanced and accurate in the modelling of the coupled systems of shallow water equation and solute transport model with friction term.